\documentclass[11pt]{amsart}
\usepackage{amsfonts}
\usepackage{latexsym}
\usepackage{amscd}
\usepackage{amsmath}
\usepackage{amssymb}
\usepackage{graphicx}
\usepackage[mathscr]{eucal}
\usepackage{xypic}
\usepackage{color}

\newtheorem{teo}{Theorem}[section]
\newtheorem{coro}[teo]{Corollary}
\newtheorem{lema}[teo]{Lemma}

\newtheorem{ejem}[teo]{Example}

\newtheorem{question}[teo]{Question}

\def\cantor{2^{\nat}}

\def\nat{\mathbb{N}}

\def\nin{\not\in}

\def\su{\subseteq}

\def\w1{\omega_1}

\begin{document}

\title{Cardinality of the Ellis semigroup on compact metric countable spaces}

\author{S. Garc\'{\i}a-Ferreira}
\address{Centro de Ciencias Matem\'aticas, Universidad Nacional
Aut\'onoma de M\'exico,  Campus Morelia,  Morelia 58190, Michoac\'an, M\'exico}
\email{sgarcia@matmor.unam.mx}

\author{Y. Rodr\'iguez-L\'opez}
\address{Secci\'on de Matem\'aticas, Universidad Nacional Experimental Polit\'ecnica ``Antonio Jose de Sucre", Barquisimeto, Venezuela}
\email{yrodriguez@unexpo.edu.ve}

\author{C. Uzc\'ategui}
\address{Escuela de Matem\'aticas, Facultad de Ciencias, Universidad Industrial de
Santander, Ciudad Universitaria, Carrera 27 Calle 9, Bucaramanga,
Santander, A.A. 678, Colombia and
Departamento de Matem\'{a}ticas,
Facultad de Ciencias, Universidad de los Andes, M\'erida 5101,
Venezuela}
\email{ cuzcatea@saber.uis.edu.co}

\thanks{The research of the first listed author was supported by  CONACYT grant no. 176202 and PAPIIT
grant no. IN-101911. The research of the second listed author was
supported by Universidad Nacional Experimental Polit\'ecnica
``Antonio Jose de Sucre", Venezuela. The third author thanks the financial assistance of VIE of Universidad Industrial de Santander. Hospitality and financial
support received from
the  {\em Centro de Ciencias Matem\'aticas de la
Universidad Aut\'onoma de M\'exico} are gratefully acknowledged}



\dedicatory{}


\begin{abstract} Let $E(X,f)$ be the Ellis semigroup of a  dynamical system   $(X,f)$  where  $X$ is  a compact   metric space.
We analyze the cardinality of $E(X,f)$ for  a   compact  countable metric space $X$.  A characterization  when
  $E(X,f)$ and  $E(X,f)^* = E(X,f) \setminus \{ f^n : n \in \mathbb{N}\}$ are both finite is given.
We show that if the collection  of all periods of the periodic points of $(X,f)$ is infinite, then
$E(X,f)$ has size  $2^{\aleph_0}$. It is also proved that if $(X,f)$  has a point with a  dense orbit and all elements of
$E(X,f)$ are continuous, then $|E(X,f)| \leq |X|$. For  dynamical systems of the form $(\omega^2 +1,f)$, we show that if there is a point with a dense orbit,  then
all  elements of   $E(\omega^2+1,f)$  are continuous functions. We present several examples of dynamical systems which have a point with a dense orbit. Such systems provide  examples where  $E(\omega^2+1,f)$ and $\omega^2+1$ are homeomorphic but not algebraically homeomorphic, where $\omega^2+1$ is taken with the usual ordinal addition  as semigroup operation.
\end{abstract}

\maketitle

\section{Introduction}

We start the paper by fixing some standard notions and terminology. Our
dynamical systems $(X,f)$'s will consist of a compact metric space  $X$
 and a continuous function $f:X\rightarrow X$.
The {\it orbit} of $x$, denoted by $\mathcal O_f(x)$, is
the set $\{f^n(x):n\in\mathbb N\}$, where $f^n$ is $f$ composed
with itself $n$-times. A point $x\in X$ is called a
{\it periodic point} of $f$ if there exists $n\geq 1$ such that
$f^n(x)=x$, and its {\it period} is $s=\min\{n\in\mathbb N:f^n(x)=x\}$.
Let $P_f$ be the set
of all periods of the periodic points of $(X,f)$.
A point $x$ is called {\it eventually periodic} if its orbit is finite.
The $\omega-${\it limit set} of $x\in X$, denoted
by $\omega_f(x)$, is the set of points $y\in X$ for which there
exists an increasing sequence $(n_k)_{k\in\mathbb N}$ such that
$f^{n_k}(x)\rightarrow y$. Observe that for each $y\in\mathcal O_f(x)$, we have that
$\omega_f(y)=\omega_f(x)$. If $\mathcal{O}_f(y)$ contains a
periodic point $x$, then $\omega_f(y)=\mathcal{O}_f(x)$. We denote
by $\mathcal{N}(x)$ the collection of all  neighborhoods of
$x\in X$.
The set of all accumulation points of $X$ will be
denoted by $X'$. For a successor ordinal $\alpha=\beta +1$,
$X^{(\alpha)}=(X^{(\beta)})'$ and for a limit ordinal $\alpha$,
$X^{(\alpha)}=\bigcap_{\beta<\alpha}X^{(\beta)}$.
The {\it Cantor-Bendixson rank of $x\in X$}, denoted by $CB(x)$, is the first ordinal $\alpha<\omega_1$
such that $x\in X^{(\alpha)}$ and $x\notin X^{\alpha+1}$. The {\it  Cantor-Bendixson rank} of $X$
is the first ordinal $\alpha<\omega_1$  for which $X^{(\alpha)}=\emptyset$.
The Stone-\v{C}ech compactification $\beta(\mathbb N)$ of
$\mathbb N$ with the discrete topology will be identified
with the set of ultrafilters over $\mathbb N$. Its remainder is denoted by
${\mathbb N}^*= \beta({\mathbb N})\setminus \mathbb{N}$ is  the
set of all free ultrafilters on $\mathbb N$, where, as usual, each
natural number $n$ is identified with  the fixed ultrafilter
consisting of all subsets of $\mathbb N$ containing $n$.
For each $A\subseteq \mathbb N$, $A^*$ denotes
the collection of all $p\in \mathbb{N}^*$ such that $A\in p$. If $A, B \in \mathcal{P}(\mathbb{N})$ and $A\setminus B$
is finite, then we will write $A\subseteq^* B$, and  $A=^*B$ when  $A\subseteq^* B$ and
$B\subseteq^* A$.

\medskip

Given a discrete dynamical system $(X,f)$ its {\it Ellis semigroup}, denoted $E(X,f)$, is
defined as  the pointwise closure of $\{f^n:\; n\in \mathbb N\}$
in the compact space $X^X$ with composition of functions as its
algebraic operation.  The Ellis semigroup is
equipped with the topology inhered from the product space $X^X.$
The Ellis semigroup of a discrete dynamical system was introduced by R. Ellis in \cite{ell} and has been very useful to study the topological behavior of the
dynamical systems.  The article \cite{glasner} offers an excellent survey concerning applications of the Ellis semigroup. In the paper \cite{GarciaSanchis}, the authors initiated the study of  the continuity of the elements of $E(X,f)^*$. For instance, they point out that  if $X$ is a convergent sequence with its limit point, then all the elements of $E(X,f)$ are either continuous or discontinuous. On the other hand, P. Szuca \cite{Szuca} showed that if $X = [0,1]$ and $f^p$ is continuous for some $p \in \mathbb{N}^*$, then all the elements of $E([0,1],f)$ are continuous. In the same direction,  some results are obtained in \cite{gf}  when the phase space is the Cantor set.
The main  tool that have been used in all these investigations   is the combinatorial properties of the ultrafilters on $\mathbb{N}$. Certainly,
the Ellis semigroup can be described  in terms of
the notion of $p-$limits where $p$  is an ultrafilter on the natural number $\mathbb N$.
Indeed, given $p\in \mathbb{N}^*$ and a sequence $(x_n)_{n\in\mathbb{N}}$
in a space $X$, we say that a point $x\in X$ is the $p-${\it limit
point} of the sequence, in symbols $x =p-\lim_{n\rightarrow
\infty}x_n$, if for every neighborhood $V$ of $x$,
$\{n\in\mathbb{N}: f^n(x)\in V\} \in p$.
Observe that a point $x\in X$ is an accumulation point of a
countable set $\{x_n:\,n\in\mathbb{N}\}$ of $X$ iff there is $p\in
\mathbb{N}^*$ such that $x = p-\lim_{n\rightarrow \infty}x_n$.

\medskip

The notion of a
$p-$limit point has been used in several branches of mathematics (see for
instance \cite{Be} and \cite[p. 179]{fu}).
A. Blass \cite{Bla} and N. Hindman \cite{hi} formally established
the connection between  ``the iteration in to\-po\-lo\-gi\-cal dynamics'' and
``the convergence with respect to an ultrafilter'' by considering a more
general iteration of the function $f$ as follows: Let $X$ be
 a compact  space and $f : X
\rightarrow  X$ a continuous function. For $p\in\mathbb{N}^*$, the
$p-$iterate of $f$ is the function $f^p: X\rightarrow X$ defined
by \[ f^p(x) = p-\lim_{n\rightarrow \infty} f^n(x),
\]
for all $x\in X$. The description of the Ellis semigroup and its operation in
terms of the $p-$iterates are the following:
\[
\begin{array}{rcl}
E(X,f) & = &\{f^p: p\in \beta\mathbb N \}\\
\\
f^p\circ f^q & =& f^{q+p}\;\; \mbox{for each $p,q\in\beta\mathbb
N$ (see \cite{Bla}, \cite{hi})}.
\end{array}
\]
We will use the following notation
\[
E(X,f)^* :=E(X,f)\setminus\{f^n:n\in\mathbb N\}.
\]
Besides, we have that $\omega_f(x)=\{f^p(x):p\in\mathbb N^*\}$ for each $x \in X$.

\medskip

In this paper we  are interested on the  cardinality of the Ellis semigroup  $E(X,f)$.
The work of A. K\"{o}hler \cite{Ko} and M. E. Glasner and Megrehisvili \cite{GM06} contain very interesting
results about when $E(X,f)$ has cardinality at most $2^{\aleph_0}$ (the so called {\em tame} dynamical systems).  In \cite{GM06} it is established the
Bourgain-Fremlin-Talagrand dichotomy for dynamical systems:  either $|E(X,f)|\leq 2^{\aleph_0} $  or $E(X,f)$ contains a copy of $\beta\nat$ and  $|E(X,f)|= 2^{2^{\aleph_0}}$.
We will be mostly concerned with the case when $X$ is countable compact metrizable space. In this case, as the cardinality of $E(X,f)$  is obviously   bounded for the cardinality of $X^X$, then $|E(X,f)|\leq 2^{\aleph_0}$.  Moreover, as  $E(X,f)$ is a separable metric space, the classical perfect set theorem says that $E(X,f)$ is either (at most) countable or has cardinality  $2^{\aleph_0} $.  Thus a natural question is to determine conditions under which each of those alternatives hold.  For instance, we characterize when  $E(X,f)$ and  $E(X,f)^*$ are finite.  We prove that if
the set of all periods of the periodic points of $(X,f)$ is infinite, then
$|E(X,f)|= 2^{\aleph_0}$. Concerning upper bounds, we prove that  if $(X,f)$ has a dense orbit and
all elements of the Ellis semigroup are continuous, then $|E(X,f)^*| \leq |X|$.
In the third section, it is shown that when $(\omega^2+1,f)$ has a point with a dense orbit, then
$E(\omega^2+1,f)$  is countable and contains only continuous functions.
An example of a continuous function $f: \omega^3+1 \to \omega^3+1$ such that
$E(\omega^3 +1,f)^*$ contains only discontinuous functions is given in section fourth.
Additionally, we provide  examples of dynamical systems with a dense orbit when the phase space is $\omega^2+1$. These dynamical  systems  also illustrate that  $E(\omega^2+1,f)$ and $\omega^2+1$ can be  homeomorphic but not algebraically homeomorphic (where $\omega^2+1$ is equipped with the standard ordinal addition as semigroup operation). 

\section{Cardinality of the Ellis semigroup}

To start this section we state several  auxiliary results that were proved in  \cite{gru}.

\begin{lema}\label{piteradaperiodico2} Let $(X,f)$ be a dynamical system and $x \in X$.
\begin{itemize}
\item[(i)] Assume that $x$ is periodic with period $n$ and let $l < n$. Then,   $p \in \big(n
\mathbb{N} + l\big)^*$ iff  $f^p(x) = f^l(x)$.

\item[(ii)]  Suppose that  $x$  is eventually periodic and that $m \in
\nat$  is the smallest positive integer such that $f^m(x)$ is a
periodic point. If $n$ is the period of $f^m(x)$ and $p \in
\big(n \mathbb{N} + l\big)^*$ for some $l < n$, then $f^p(x) =
f^l(f^{nj}(x))$ where $j=\min\{i:\, m\leq ni+l\}$.

\item[(iii)]  Suppose that the orbit of  $x$ is  infinite and $\omega_f(x)=\mathcal O_f(y)$
for some periodic point $y \in X$ with period $n$. If  $p,\,\,q\in(n\mathbb N+l)^*$ for some  $l<n$, then  $f^p(x)=f^q(x)$.

\item[(iv)]   $f^p(f^n(x))=f^n(f^p(x))$  for every  $n\in\nat$, $x\in X$ and every $p\in\nat^*$.
\end{itemize}
\end{lema}

The next theorem provides a necessary and sufficient conditions for $E(X,f)$ to be
finite, its elementary proof is omitted.

\begin{teo} \label{Ellis finito} Let $(X,f)$ be a dynamical system. Then $E(X,f)$ is finite iff
there exists  $M>0$ such that
$|\mathcal O_f(x)|<M$ for each $x\in X$.
\end{teo}

Similarly, we will show a sufficient and necessary condition for $E(X,f)^*$
to be finite.

\begin{teo}\label{p(f)finitoasterisco} Let $(X,f)$ be a dynamical system. $E(X,f)^*$ is finite
iff there is $M \in \mathbb{N}$ such that $|\omega_f(x)|\leq M$ for each $x\in X$.
\end{teo}

\proof Assume that  $E(X,f)^*$ is finite and let $M=|E(X,f)^*|$. Since $\omega_f(x)=\{f^p(x):p\in\mathbb N^*\}$,
we must have that $|\omega_f(x)| \leq M$ for each $x\in X$.

Conversely, assume that there is $M \in \mathbb{N}$ so that
$|\omega_f(x)| \leq M$ for each $x\in X$. Then,  we know
 that every point of $\omega_f(x)$ is periodic for each $x\in X$
and so $|P_f|\leq M$. Hence, let $P_f=\{b_1, \cdots, b_n\}$ for some $n \leq M$. Define
$\phi:\mathbb N^*\rightarrow \prod_{i=1}^n\{0,\cdots, b_i-1\}$ by
$\phi(p)=(j_1,....,j_n)$ provided that   $p\in (b_i\nat+j_i)^*$ for each
$1\leq i\leq n$, for every  $p\in\mathbb N^*$. To see that $E(X,f)^*$ is
finite, it suffices to show that $\phi(p)=\phi(q)$ iff $f^p=f^q$, for $p, q \in \mathbb{N}^*$. But
this follows directly from clauses  $(i)$ and $(iii)$ of Lemma \ref{piteradaperiodico2}.
\endproof

It is noteworthy that  $E(X,f)^*$ could be finite and $E(X,f)$ could be infinite. For instance, if  $X$ is a convergent sequence
with its limit point and $f$ is the shift function, then $E(X,f)$
is infinite and $E(X,f)^*$ consists of only  one point.
The next corollary  follows directly  from Theorem \ref{p(f)finitoasterisco}.

\begin{coro} Let $(X,f)$ be a dynamical system. If $E(X,f)^*$ is finite,
then $P_f$ is finite.
\end{coro}

The converse of the previous corollary is not true,   we shall describe an example of a dynamical system $(X,f)$
such that $P_f$ is finite and $E(X,f)^*$ infinite (Examples \ref{orbita densa 2} and \ref{orbitad-densa3}).

\medskip

So far we have dealt only with conditions implying that $E(X,f)$ or $E(X,f)^*$  are finite. Now we show that  $E(X,f)$ has cardinality at least
$2^{\aleph_0}$, whenever  $P_f$ is infinite. For that end we need  the general form of the Chinese Remainder Theorem
(for more properties about this theorem see for instance \cite{ro}).

\begin{lema}\label{ChRT}
The system of equations
\[
x\equiv r_i \,(\mbox{mod} \; m_i),
\]
for $i=0, \cdots, k$, has an integer solution $x$ iff
$gcd(m_i,m_j)$ divides $r_i-r_j$ for all $i\neq j$.
\end{lema}

The following lemma  is very important for the estimation of the cardinality of the Ellis semigroup by using ultrafilters on $\mathbb{N}$.

\begin{lema}
\label{extension}
Let   $k\geq 1$ and $m_1, \cdots, m_k$ be  positive integers, $0\leq r_i <m_i$. Suppose that the following system $E$  of equations has a solution:
\begin{equation}
\label{ecu1}
x\equiv r_i\,(\mbox{mod} \; m_i)\;\;\mbox{for $i=1, \cdots, k$.}
\end{equation}
Then for every  infinite $A\su \nat $ there are an infinite subset $B$ of $A$ and positive integers $s_1<s_2$ such that the equation system $E\cup\{x\equiv s_i (\mbox{mod}\; m) : i=1,2 \}$ has a solution for all $m \in B$ with $m > s_2$.
\end{lema}

\proof By the Pigeon Hole Principle, there is an infinite $B\su A$  such that if $m, m'\in B$, then $gcd(m,m_i)=gcd(m',m_i)$ for all $1\leq i\leq k$.
Observe that  the number  $l_i= gcd(m,m_i)$ does not depend on the choice of $m\in B$, for each $1\leq i\leq k$. Consider the following system of equations:
\begin{equation}
\label{ecu2}
x\equiv r_i\,(\mbox{mod} \; l_i)\;\;\mbox{for $i=1, \cdots, k$.}
\end{equation}
Since the system of equations \eqref{ecu1} has a solution,  by Lemma \ref{ChRT}, we have that $gcd(m_i,m_j)$  divides $r_i-r_j$ for every $1\leq i\leq k$. As $gcd(l_i,l_j)$ divides $gcd(m_i,m_j)$, then by Lemma \ref{ChRT} the system \eqref{ecu2} has a solution. Choose two such solutions $s_1< s_2$. If $m \in B$ satisfies that $s_2<m$,
again by Lemma \ref{ChRT}, then the system $E\cup\{x\equiv s_i (\mbox{mod}\; m) : i=1,2 \}$ has a solution.
\endproof

Next, we will see how the machinery of Number Theory  and the ultrafilters on $\mathbb{N}$ works.

\begin{teo}\label{infipf}
Let $(X,f)$ be a dynamical system. If $P_f$ is infinite, then
 $E(X,f)$ has no isolated points. In particular, if $X$ is countable and  $P_f$ is infinite, then $E(X,f)$ is homeomorphic to  $\cantor$.
\end{teo}

\proof  First, we shall prove that $E(X,f)$ has no isolated points.
Let $p\in \beta\nat$  and $V$ be an open set  in $X^X$ such that $f^p\in V = \{g\in X^X:\; g(x_i)\in V_i\; \mbox{for $i=1,\cdots , k$}\}$ where    $x_1, \cdots, x_k\in X$  and  $V_i \cdots, V_k$ are open subsets of $X$. It suffices to show that there is $q\in \beta\nat$ such that $f^q\neq f^p$ and $f^q\in V$. We need to consider two cases.

Case 1:  There is $1 \leq i \leq k$  such that $x_i$ has infinite orbit.  If $B_j=\{n\in\nat:\, f^n(x_j)\in V_j\}$, then $p\in B_j^*$ for every $j\leq k$. Hence, we can find  $n\in B_0\cap B_1\cap\cdots\cap B_k$ so that $f^n(x_i)\neq f^p(x_i)$. Notice that  $f^n \in V$.

\medskip

Case 2:  $x_i$ has finite orbit for every $1 \leq i \leq k$.  So each $x_i$ must be  eventually periodic of some period $m_i$. Pick  $r_i<m_i$ such that $p\in (m_i\nat+r_i)^*$  for each $1 \leq i\leq k$.
By Lemma \ref{piteradaperiodico2} (ii), $f^q(x_i)=f^p(x_i)$ whenever $ q\in (m_i\nat+r_i)^*$ for each $1 \leq i\leq k$.
By Lemma \ref{extension}, there are $m\in P_f$  and $s_1<s_2<m$ such that  the system of equations $ x\equiv r_i\,(\mbox{mod} \; m_i)$ for $i=1, \cdots, k$  together with the equation $ x\equiv s_j\,(\mbox{mod} \; m)$ has a solution for $j=1,2$.  Pick $p_j\in \bigcap_{i=1}^k (m_i\nat+r_i)^* \cap (m\nat+s_j)^*$ for $j=1,2$. Let $y$ be a periodic point with period $m$.  Then $f^{p_1}(x_i)=f^{p_2}( x_i)$, for every $1 \leq i \leq k$,  and $f^{p_j}(y) = f^{s_1}(y) \neq f^{s_2}(y) = f^{p_2}(y)$. Therefore,  $f^{p_1}, f^{p_2}\in V$ and $f^{p_1}\neq f^{p_2}$.

\medskip

If $X$ is countable, then   $E(X,f)$ is homeomorphic to $\cantor$ since it is compact metrizable without isolated points and  zero dimensional
\endproof

Concerning the previous theorem, it would be also interesting to see whether the algebraic structure of $E(X,f)$ is unique when $X$ is countable and $P_f$ is infinite. We also wonder about the existence  of a dynamical system $(X,f)$ for which $P_f$ is finite and $E(X,f)$ is uncountable.

\begin{coro}\label{cantorset} For every  countable ordinal $\alpha \geq 1$, there is  a continuous function $f: \omega^\alpha+1 \rightarrow \omega^\alpha+1$ such that  $E(\omega^\alpha+1,f)$ is homeomorphic to the Cantor space $\cantor$.
\end{coro}

\proof According to Theorem \ref{infipf}, it suffices to  define a continuous function $f:  \omega^\alpha+1 \rightarrow \omega^\alpha+1$ for which  $P_f$ is infinite.  For the case when  $\alpha=1$, it is very  easy to define such a function. Suppose that  $\alpha>1$. Choose a subspace $X$  of  $\omega^\alpha+1$ so that $\omega^\alpha+1=(\omega+1) \oplus X$. If  $f$ is the function defined for the case $\alpha=1$, then we consider the function $g =f\oplus Id: (\omega+1) \oplus X \rightarrow (\omega+1) \oplus X$. Observe that $g^p=f^p\oplus Id$, for every $p\in \beta(\nat)$, which implies that  $E(\omega^\alpha +1,g)\approx E(\omega+1, f) \approx \cantor$.
\endproof

\medskip

In the next theorem, we shall bound the cardinality of the Ellis semigroup for  certain dynamical systems.

\begin{teo}\label{Ellis acotado} Let $(X,f)$ be a dynamical system such that there is  $w\in X$ with  a dense orbit.
Suppose that  $f^p$ is continuous for
every $p\in \nat^*$. Then  $f^p=f^q$ iff $f^p(w)=f^q(w)$, for every $p, q \in \mathbb{N}^*$.
In particular, $|E(X,f)^*| \leq |X|$.
\end{teo}

\proof Let $p\in\mathbb N^*$. It suffices to prove that $f^p$ is completely determined by its value at $w$. Fix  $x\in X$. First, suppose that   $x$ is an isolated point. Then, there is $n\in\mathbb N$ (depending on $x$)  such that $f^{n}(w)=x$. Thus,
$f^p(x)=f^{p}(f^{n}(w)) = f^{n}(f^{p}(w))$ for every  $p\in\nat^*$. Now, assume  otherwise that $x$ is a limit point.
Choose a sequence $(f^{m_k}(w))_{k\in \nat}$ converging to $x$. Since $f^p$ is continuous, we have that
\[
f^p(x)=\lim_{k \to \infty} f^p(f^{m_k}(w))=\lim_{k \to \infty} f^{m_k}(f^p(w)).
\]
Therefore,  $f^p$  is completely determined by $f^p(w)$.
\endproof

\section{Compact Metric Countable Spaces}

We remind the reader the classical  result that  every  compact metric countable
space  is homeomorphic to a countable ordinal with
the order topology (see \cite{ms}). In what follows, our  dynamical systems will have
$\omega^\alpha+1$ as a phase space, where $\alpha$ is a countable ordinal with $\alpha\geq 1$. For our convenience,
 $d$ will stand for the unique point of $\omega^\alpha+1$ with $CB$-rank equal to $\alpha$.

\medskip

In the next, lemma we  list some basic properties of the dynamical systems of the form  $(\omega^\alpha+1,f)$ with a dense orbit.

\begin{lema}\label{omega-alpha}
Let $(\omega^\alpha+1,f)$ be a dynamical system with $\alpha\geq 1$,
such that there exists $w \in \omega^\alpha+1$ with a dense orbit. Then, the following conditions hold:
\begin{itemize}
\item[$(i)$]   $f(y)$ is a limit point for every   $y \in (\omega^\alpha+1)'$.

\item[$(ii)$]  $w$ is isolated and its orbit  consists of all isolated points of $\omega^\alpha+1$.

\item[$(iii)$]  The range of $f$  is $\omega^\alpha+1\setminus\{w\}$.

\item [$(iv)$] If $x \in (\omega^\alpha+1)'$, then $\emptyset \neq f^{-1}(x) \subseteq (\omega^\alpha+1)'$.

\item [$(v)$] If $CB(z)=1$, then $z$ is not periodic.
\end{itemize}
\end{lema}

\proof $(i)$. Let  $y \in (\omega^\alpha+1)'$ and suppose  $f(y)$ is  isolated to get a contradiction.  Since $y$ is a limit point, then $\{z\in \omega^\alpha+1:
\; f(z)=f(y)\}$ is open and infinite. As the orbit of $w$ is dense, there are $k<l$ such that $f(f^{k}(w))=f(f^{l}(w))=f(y)$, thus the orbit of $w$ is finite, which is a contradiction.

$(ii)$.  Since the orbit of $w$ is dense, then it contains all isolated points of $\omega^\alpha+1$ and  by $(i)$ no limit point belongs to the orbit of $w$.

$(iii)$. Since the orbit of $w$ is dense, then it is clear that every limit point belongs to the range of $f$ since it is compact.  If  $f(y)=w$ for some $y\in X$, then $y$ must be isolated and so $w$ is  periodic, which is impossible. So, $w$ is not in the range of $f$.

$(iv)$. Let $x \in (\omega^\alpha+1)'$. By clause $(iii)$,  we have that $f^{-1}(x) \neq \emptyset$. Suppose that $y\in f^{-1}(x)$ is isolated.
Since the orbit of $w$ is dense,  there is $n\in\nat$ such that $f^n(w)=y$ and thus  $f^{n+1}(w)=x$.   By $(i)$,  $f^m(w)$ is a limit point  for all $m\geq n+1$, but this contradicts the fact that the orbit of $w$ is dense.

$(v)$ Suppose that  $z$ is a periodic point of period $l$.
For each $0\leq j<l$, fix disjoint  clopen sets $V_j$ such that $V_j\cap \mathcal{O}_f(z)=\{f^j(z)\}$ and also $z$ is the only limit point in $V_0$ (since $CB(z)=1$).
Moreover, we can also assume that $f[V_{l-1}]\subseteq V_0$ and $f[V_i]\subseteq V_{i+1}$ for $0<i<l-1$ (notice that we cannot ask that $f[V_0]\subseteq V_1$).
Let  $V=V_0\cup\cdots \cup V_{l-1}$ .
Since the orbit of $w$ is  dense, then there is $j<l$ such that
$A=\{n\in\mathbb N: f^{n}(w) \in V_j \text{ and }f^{n+1}(w)\notin V \}$ is infinite.
From the assumptions about the $V_i$'s, it is clear that $j=0$. Let $\{n_k : k\in\mathbb N\}$ be an infinite subset of $A$ such that  $f^{n_k}(w)\rightarrow z$.
However, $f (w)\in V$ but $f(f^{n_k}(w))\notin V$ for all $k$, contradicting the continuity of $f$.
Consequently, $z$ is not periodic.
\endproof

\begin{lema}\label{omega-alpha2} Let $(\omega^\alpha+1,f)$ be a dynamical system with $\alpha\geq 1$
such that there exists $w\in \omega^\alpha+1$ with a dense orbit.

\begin{itemize}
\item[$(i)$]   Suppose  $x\in(\omega^\alpha+1)'$ is such that   $CB(y)<CB(x)<\alpha$ for every  $y\in f^{-1}(x)$.
If $(x_n)_{ n \in \nat }$ is a sequence  such that  $ x_n \to x$ , $CB(x_n)<CB(x)$, for each $n$, and $\sup\{CB(x_n) : n \in \nat \} = CB(x)$, then there is $N \in \nat$ such that $CB(z)<CB(x_n)$ for every  $z\in f^{-1}(x_n)$  and $n \geq N$.

\item[$(ii)$]  $f(d)=d$.
\end{itemize}
\end{lema}

\proof
$(i)$ Let $(x_n)_{ n \in \nat }$ be a sequence  such that $ x_n \to x$,  $CB(x_n)<CB(x)$ for all $n$  and $CB(x_n)$ converges to $CB(x)$.
We proceed by contradiction. By  Lemma \ref{omega-alpha}, and passing to a subsequence if necessary, there is  a sequence  $(z_n)_{n\in\mathbb N}$ so that $z_n\in f^{-1}(x_n)$
and $CB(z_n)\geq CB(x_n)$  for each $n\in\mathbb N$. Choose a sequence $(n_k)_{k\in\mathbb N}$ and
$y$ such that  $z_{n_k}\rightarrow y$. Clearly $y\in f^{-1}(x)$. Since $CB(z_{n_k})\geq CB(x_{n_k})$, then  $CB(y)\geq CB(x)$, which is a contradiction.

$(ii)$ Assume that $f(d)=y$ for some $y\in(\omega^\alpha+1)'\setminus \{d\}$ with $CB(y)<\alpha$.   Let $d_n \to d$ such that $CB(d_n)<\alpha$ for all $n$ and $\sup_n CB(d_n)=\alpha$.
By the continuity of $f$ at $d$, we can assume without loss of generality that $f(d_n)\neq d$ for all $n$.
We claim that there is $n \in\mathbb N$ such that $CB(z)<CB(d_n)$ for all $z\in f^{-1}(d_n)$.
Otherwise, for each $n\in \mathbb N$, there is $z_n\in f^{-1}(d_n)$ such that $CB(z_n)\geq CB(d_n)$. Therefore, passing to a subsequence if it is necessary,  $z_n\to d$ and $f(z_n)\to d$, then $f(d)=d$ which  contradicts our assumption.

Fix $n_0\in \mathbb N$  such that $CB(z)<CB(d_{n_0})$ for all  $z\in f^{-1}(d_{n_0})$.  Let $y_1=d_{n_0}$. We can apply part $(i)$ and obtain a point $y_2$ such that
$CB(y_2)<CB(y_1)$ and $CB(z) < CB(y_2)$ for every  $z\in f^{-1}(y_2)$. This process has to end in a finite number of steps  when we reach a point $y_k$ such that   $CB(z) < CB(y_k) = 1$ for every  $z\in f^{-1}(y_k)$. But this contradicts part  $(iv)$ of Lemma \ref{omega-alpha}. Therefore, $f(d) = d$.
\endproof

\begin{teo}\label{orbita-densa-omega2} Let $(\omega^2 +1,f)$ be a dynamical system
such that there exists $w\in \omega^2+1$ with a dense orbit.
Then $f^p$ is continuous, for every $p\in \nat^*$, and  $E(\omega^2+1,f)$ is  homeomorphic to $\omega^2+1$.
\end{teo}

\proof  Fix $p\in\nat^*$.  Let $\{ d_n : n \in \mathbb{N} \}$ the set of all points of $X$ with $CB$-rank equal to $1$. According to Lemma \ref{omega-alpha2}, we know that $f^p(d) = d$ and,
by Lemma \ref{omega-alpha}, we have that $f(d_n)$ is a limit point for each $n \in \nat$. Thus, we obtain that $O_f(d_ n)\subseteq \{d_m:\; m\in\nat\}\cup \{d\}$ for each $n\in \nat$. Besides, by Lemma \ref{omega-alpha} $(v)$, we know that $d_n$ cannot be  periodic for all $n \in \nat$.

First, we show that  $f^p(d_n) = d$ for each $n \in \nat$. Fix $n \in \nat$. To have this done we consider two cases: $(a)$ Suppose that $\mathcal{O}_f(d_n)$ is finite. Then  $d_n$ is eventually periodic. Since $d_m$ is not periodic for every $m \in \nat$, then $d$ belongs to the orbit of $d_n$ and thus  $f^p(d_n)=d$.
$(b)$ Suppose that $\mathcal{O}_f(d_n)$ is infinite. Since   $O_f(d_ n)\subseteq \{d_m:\; m\in\nat\}\cup \{d\}$,  then $f^m(d_n) \rightarrow d$ and   thus  $f^p(d_n)=d$.

We are ready to prove the continuity of  $f^p$. Let $y\neq d$ be a limit point and $(y_n)_{n \in \nat}$ be a
sequence of isolated points converging to $y$.  We have already shown that $f^p(y)=d$. For each $n \in \nat$, put $y_n=f^{k_n}(w)$ for some $k_n \in \nat$. Hence, we have that
$f^p(y_n)=f^{k_n}(f^p(w))$ for every $n \in \nat$. Since  $f^p(w)$ is a limit point, then $f^p(y_n)\rightarrow
d$.   Finally, if $y_n\rightarrow d$, then regardless of whether each $y_n$ is isolated or not, we have that $f^p(y_n)\rightarrow d$.

Now we will verify that  $E(\omega^2 +1,f)$  is homeomorphic to $\omega^2 +1$. For each $n \in \nat$, let us fix a clopen subset $V_n$ of $\omega^2 +1$ such that $d_n$ is the only limit point in $V_n$. Consider the following sets:
\[
B_n=\{i \in  \nat :\; f^i(w)\in V_n \}.
\]
We claim that:

$(i)$  If $p,q\in B_n ^*$  for all $n \in \nat$, then $f^p=f^q$ and $f^p(w) = d_n$.

$(ii)$ If $p,q\nin B_n^*$ for all $n \in \nat$, then $f^p=f^q =$ the constant function with constant value equal to $d$.

\medskip

For every $n \in \nat$ fix $p_n \in B_n^*$ and fix $q \notin \bigcup_{n \in \nat}B_n^* $.  Then we define $F: \omega^2+1\rightarrow E(X,f)$ as follows:
\[
F(x) = \left\{\begin{array}{lcl} f^m & & \mbox{if $x=f^m(w)$}
\\
\\
f^{p_n} & & \mbox{if $x=d_n$ \ for some} \ n \in \nat \\
\\
f^q & & \mbox{if $x=d$.}
\end{array}
\right.
\]
It is not hard to prove that  $F$ is an homeomorphism.
\endproof

\medskip

Theorem  \ref{orbita-densa-omega2} is false for the space $\omega^3 +1$. Certainly, in Example \ref{orbitad-densa3}, we shall construct a continuous function $f$ on $\omega^3 +1$ with a dense orbit
for which  $f^p$ is discontinuous on $\omega^3 +1$, for all $p\in\nat^*$.

\medskip

We end this section by making some remarks about the algebraic structure of  $E(\omega^2+1,f)$ when it has a dense orbit.  It follows from the proof of the previous theorem that the semigroup operation of  $E(\omega^2+1,f)$ satisfies the following properties:

 \begin{enumerate}
 \item if $p, q \in \nat^*$, then $f^p \circ f^q =$ the constant function with constant value equal to $d$.

 \item If $p \in B_n^*$, for some $n \in \nat$, and $k \in \nat$, then $f^p \circ f^k(d_m) = f^k(f^p(d_m)) = f^k(d) = d$ for every $m \in \nat$, and $f^p \circ f^k(w) = f^k(d_n) = f^k(d_n)$.

 \item If $p \notin \bigcup_{n \in \nat}B_n^*$ and $k \in \nat$, then $f^p \circ f^k =$ the constant function with constant value  $d$.
\end{enumerate}

On the other hand, we observe that each ordinal of the form $\omega^\alpha$ is a semigroup under the usual ordinal addition. We can extend this semigroup operation to $\omega^\alpha+1$ by simply declaring $\omega^\alpha+\beta=\beta+\omega^\alpha=\omega^\alpha$ for all $\beta \in \omega^\alpha+ 1$. Under this operation, $\omega^\alpha+ 1$ is a semigroup and the operation is left continuous.  It is known that $E(X,f)$ is  also a left continuous semigroup for any dynamical system $(X,f)$.  Considering a dynamical system $(\omega^2+1,f)$ satisfying the conditions of Theorem \ref{orbita-densa-omega2}, we can see that $E(\omega^2+1,f)$ is not topological isomorphic (i. e., both isomorphic and homeomorphic)  to $\omega^2 +1$ equipped with the above operation: Indeed, let $\{ d_n : n \in \nat \}$  be the set of  limit points of $\omega^2+1$ with $CB$-rank equal to $1$. We have that, under ordinal addition, $d_n+d_n< d$, for every $n \in \nat$, but $f^{p}\circ f^{q}=$ constant function with value $d$ provided that $p, q \in \nat^*.$

\section{Examples }

In this final section,  we shall present some examples and counterexamples  of dynamical systems which contain a dense orbit and have very interesting topological properties. The phase spaces of these dynamical systems will be $\omega^2+1$ and $\omega^3+1$. The construction of the first example is based in the following combinatorial lemma.

\begin{lema}\label{una-sola-orbita}
There is a family $\{A_k:k\in\mathbb N\}$ of pairwise disjoint
infinite subsets of $\mathbb N$ and an injective function
$h:\bigcup_{k\in \nat} A_k\rightarrow \bigcup_{k \in \nat} A_k$ such that:

\begin{enumerate}
\item For each $F\in[A_{2k+2}]^{<\omega}$ there is $E\in[A_{2k}]^{<\omega}$
such that $h[A_{2k}\setminus E] \subseteq A_{2k+2}\setminus F$.

\item For each $F\in[A_{2k+1}]^{<\omega}$ there is
$E\in[A_{2k+3}]^{<\omega}$ such that $h[A_{2k+3}\setminus E]\subseteq
A_{2k+1}\setminus F$.

\item For each $E\in [A_0]^{<\omega}$
there are  $H\in [A_1]^{<\omega}$ such that
 $h[A_1\setminus H]\subseteq A_0\setminus E$.

\item There exists $x_0\in A_0$ such that
$\mathcal{O}_h(x_0)=\bigcup_{k\in \nat} A_k$.
\end{enumerate}
\end{lema}

\proof To define the function $h$ and the sets $A_k$'s we need the following:

\smallskip

Let $\{d^m_n: n,m\in \mathbb{N} \}$ be a faithfully  enumeration of
$\mathbb{N}$. Thus we have that  $\{d^m_n: n\in \mathbb{N}\}$ is infinite
for each $n\in \mathbb{N}$. Define $a_0=0$ and $a_{n+1} = a_n +
2(n+1)+1$ for each $n\in\mathbb N$.

\smallskip

Our point will be $x_0 = d^0_0$. Then our main task is to describe the orbit
of the point $d^0_0$ under $h$, which is defined by
the following  rule:

\medskip

\noindent$d^0_0\rightarrow d^2_1\rightarrow d^1_2\rightarrow d^0_3$\\
\\
$d^0_3\rightarrow d^2_4\rightarrow d^4_5\rightarrow d^3_6\rightarrow d^1_7\rightarrow d^0_8$\\
\\
$d^0_8\rightarrow d^2_9\rightarrow d^4_{10}\rightarrow d^6_{11}\rightarrow d^5_{12}\rightarrow d^3_{13}\rightarrow d^1_{14}\rightarrow d^0_{15}$\\
\\
$d^0_{15}\rightarrow d^2_{16}\rightarrow d^4_{17}\rightarrow
d^6_{18}\rightarrow d^8_{19}\rightarrow d^7_{20}\rightarrow
d^5_{21}\rightarrow d^3_{22}\rightarrow d^1_{23}\rightarrow
d^0_{24}$\\
$$\vdots$$\\
$d^0_{a_n}\rightarrow\cdots d^{2k}_{a_n+k}\rightarrow d^{2(k+1)}_{a_n+k+1}\cdots
d^{2(n+1)}_{a_n+n+1}\rightarrow d^{2n+1}_{a_n+n+2}\cdots
d^{2k+1}_{a_n+2(n+1)-k} \to d^{2(k-1)+1}_{a_n+2(n+1)-(k-1)}\cdots d^0_{a_{n+1}}$ \\
$$\vdots$$

\medskip

\noindent More precisely,  $h$ is defined  as follows:
\begin{enumerate}
  \item[$(i)$]   For each $k\in\mathbb N$, $n>k-1$, $h(d^{2k}_{a_n+k})=d^{2(k+1)}_{a_n+(k+1)}$.
  \item[$(ii)$]  For each $n\in \mathbb N$, $h(d^{2(n+1)}_{a_n+(n+1)})=d^{2n+ 1}_{a_n+ n+2}$.
  \item[$(iii)$] For each $k>0$, $n\geq k$, $h(d^{2k+1}_{a_n+2(n+1)-k})=d^{2(k-1)+1}_{a_n+2(n+1)-(k-1)}$.
  \item[$(iv)$]  For each $n\in\mathbb N$, $h(d^{1}_{a_n+2(n+1)})=d^0_{a_{n+1}}$.
  \item[$(v)$]   For each $n\in\nat$, $h^{a_n+n+1}(d^0_0)=d^{2n+1}_{a_n+n+2}$,  $h^{a_n+k}(d^0_0)=d^{2k}_{a_n+k}$ for all $0\leq k<n+1$, and
              $h^{a_n+k}(d^0_0)=d^{2(a_{n+1}-(a_n+k))-1}_{a_{n+1}-(a_{n+1}-(a_n+k))}$ for all $n+1<k<a_{n+1}-a_n$.
\end{enumerate}

\noindent Define $A_{0}=\{d^{0}_{a_n}: n \in\mathbb N\}$ and, for each  positive $k \in \nat$, we define
$A_{2k}=\{d^{2k}_{a_n+k}: n \geq k-1\}$ and
$A_{2k+1}=\{d^{2k+1}_{a_n+2(n+1)-k}: n\geq k\}$. It is clear that
the sets $A_k's$ are pairwise disjoint, and also they satisfy the following:

\begin{enumerate}
  \item[$(a)$] $h[A_{2k}]=^* A_{2(k+1)}$ for all $k\in\mathbb N$.
  \item[$(b)$] $h[A_{2k+1}]=^*A_{2k-1}$ for all $k>0$.
  \item[$(c)$] $h[A_1]=^* A_0$.
  \item[$(d)$] $h$ is injective.
  \item[$(e)$] $\mathcal O_h(d_0^0)=\bigcup_{m \in \nat} A_m$.
\end{enumerate}
Then we let the reader to  verify that $(1)$, $(2)$, $(3)$ and $(4)$ follow directly from $(a)$, $(b)$, $(c)$, $(d)$ and $(e)$.
\endproof

\medskip

To describe our examples the space $\omega^2+1$ as a subspace of $\mathbb{R}$ will be written as:
$$
\big(\bigcup_{m\in\mathbb N}(D_m \cup \{d_m\}) \big)\cup\{d\}.
$$
where the points of  $D_m=\{d_n^m:n\in\nat\}$ are isolated and
form a strictly increasing sequence converging to $d_m$, for each $m \in \nat$, and $(d_m)_{m \in \nat}$ is a strictly increasing sequence converging to
$d$.

\begin{ejem}\label{familia disjunta2} There is a continuous function $f: \omega^2+1 \to \omega^2+1$
 such that
 \begin{enumerate}
\item  $\mathcal O_f(d_m)$ is infinite for all $m \in \nat$, and

\item   $\mathcal{O}_f(d_0^0)$ is dense.
\item  $f$ is injective.

\item $E(\omega^2+1,f)$ is homeomorphic to  $\omega^2+1$.
\end{enumerate}
\end{ejem}

\noindent We will use the family $\{ A_m : m\in \mathbb{N} \}$ of pairwise disjoint
infinite subsets of $\mathbb N$  and the
function $h$ given in  Lemma \ref{una-sola-orbita}. For our convenience, we shall put
 $D_m = A_m$ for every $m \in \nat$. Thus, we define $f$ on $\bigcup_{m\in\mathbb N}D_m$ by following the function  $h$. To guarantee  the continuity of the function $f$
 the values on the non-isolated points are defined as follows:
\begin{enumerate}
\item  $f(d_{2k}) =  d_{2(k+1)}$ for each $k\in\mathbb N$.

\item $f(d_{2k+1}) =d_{2k-1}$ for each $k\in\mathbb N$.

\item $ f(d_1) =d_0$.

\item $f(d)= d$.
\end{enumerate}

The function  $f$ satisfies the following identities:

\begin{enumerate}
\item [$(a)$] $f[D_0]=D_2$.

\item [$(b)$] $f[D_{2k}]=^*D_{2(k+1)}$ for every $k>0$.

\item [$(c)$] $f[D_{2k+1}]=^*D_{2k-1}$ for every $k> 0$.

\item [$(d)$] $f[D_{1}]=^*D_0$.
\end{enumerate}

Let us check that $f$ is continuous on $\omega^2+1$. In fact, conditions
$(a)$, $(b)$ and $(1)$ guarantee  that $f$ is continuous at $d_{2k}$ for
all $k$. Similarly, $(c)$, $(d)$, $(2)$ and $(3)$ show the continuity at
$d_{2k+1}$ for all $k$. Finally, to prove the continuity of $f$ at $d$, notice that from
$(1)$ and $(2)$  that $f(d_{n_i})\rightarrow d$ for every increasing sequence $(n_i)_{i\in\mathbb N}$ of $\nat$.
On the other hand, suppose $d^{m_i}_{n_i}\rightarrow d$. Without loss of generality, we  may assume that the sequence
$(m_i)_{i\in\mathbb N}$ is strictly increasing. By conditions $(b)$ and $(c)$, we
conclude that $f(d^{m_i}_{n_i})\rightarrow d$. Therefore, $f$ is
continuous on $d$. Now, by condition $(4)$ of lemma \ref{una-sola-orbita}, we obtain that
$\mathcal O_f(d_0^0) = \bigcup_{m\in\mathbb N} D_m$ and so the orbit of $d_0^0$ is dense
in $\omega^2+1$. Moreover,  we have that $w_f(x)=\{d_m:m\in\mathbb{N}\}\cup\{d\}$, for all
$x \in\bigcup_{m\in\mathbb N} D_m$. By Theorem \ref{orbita-densa-omega2}, we obtain that
 $E(\omega^2+1,f)$ is homeomorphic to $\omega^2+1$, this shows clause $(4)$.
 $\Box$

\bigskip

Our next example is a dynamical systems $(\omega^2+1,f)$  which is a little different from the previous one.

\begin{ejem}\label{orbita densa 2} There is a continuous function $f: \omega^2+1 \to \omega^2+1$
such that
 \begin{enumerate}
\item  $\mathcal O_f(d_m)$ is finite for all $m \in \nat$, and

\item  $\mathcal{O}_f(d_0^0)$ is dense.
\item  $f$ is not injective.

\item $E(\omega^2+1,f)$ is homeomorphic to $\omega^2+1$.
\end{enumerate}
\end{ejem}

For our convenience,   suppose that $D_m=\{d_n^m:n\geq m-1\}$,  for each $m>0$.
We define the function $f:\omega^2+1\rightarrow \omega^2+1$ as follows:
\begin{enumerate}
\item[$(i)$] $f(d_m^{0}) = d_m^{m+1}$ for every $m\in\mathbb N$.

\item[$(ii)$]  $f(d^m_{n}) =  d^{m-1}_{n}$ for every $m>1$ and $n\geq m-1$.

\item[$(iii)$] $f(d^1_{m}) =d^0_{m+1}$ for every $m\in \mathbb N$.

\item[$(iv)$]  $f(d_0) =d$.

\item[$(v)$] $f(d_m) =d_{m-1}$ for every  $m>1$.

\item[$(vi)$] $f(d)= d$.
\end{enumerate}
From the definition of $f$ it is clear that
\begin{enumerate}
  \item[$(a)$] $f[D_m]=^*D_{m-1} $ for all $m>0$.

  \item[$(b)$] $f[D_0]=\{d^{m+1}_{m}:m \in \nat\}$.
\end{enumerate}
These conditions $(a)$ and $(b)$ guarantee  that $f$ is continuous at $\omega^2+1$. Applying again Theorem \ref{orbita-densa-omega2}, we conclude that
 $E(\omega^2+1,f)$ is homeomorphic to $\omega^2+1$
 $\Box$

\medskip

The  next counterexample will testify  that Theorem \ref{orbita-densa-omega2} is false for  $\omega^3+1$.
For our convenience,  the space $\omega^3+1$ as a subspace of $\mathbb{R}$ will be written as follows:

\medskip

The isolated points will be $D_{i,l}=\{d_{i,l}^k: k \geq l\}$, for each $i, l \in \nat$ with $l\geq i$ and
$$
\bigg(\bigcup_{i\in\nat, \ l \geq i}(D_{i,l} \cup \{d_{i,l}\})\bigg) \cup
 \{d_i : i \in \nat \} \cup\{d\}.
$$
For each $i \in \nat$ and $l \geq i$,  $(d_{i,l}^k)_{l \leq k}$ is a strictly increasing sequence converging to $d_{i,l}$;
$(d_{i,l})_{l\geq i}$ is a strictly increasing sequence converging to
$d_i$; and $(d_{i})_{i\in\nat}$ is a strictly increasing sequence converging to
$d$.

\begin{ejem}\label{orbitad-densa3} There is a continuous function $f: \omega^3+1 \to \omega^3+1$
 such that
 \begin{enumerate}
\item  $d_i$ is fixed for all $i\in\nat$.
\item  $\mathcal{O}_f(d_{0,0}^0)$ is dense.
\item  $f^p$ is discontinuous for all  $p\in\nat^*$.
\end{enumerate}

\end{ejem}

To define the function $f$ on the isolated points we start from the point $d^0_{0,0}$ and then we shall describe its orbit
following the rule:

\medskip

For each $ 1 \leq i \leq  k$,

\medskip

\noindent$d^k_{i,i}\rightarrow d^k_{i-1,k},$  \\
\\
$d^0_{0,0}\rightarrow d^1_{1,1}\rightarrow d^1_{0,1}\rightarrow d^1_{0,0}$\\
\\
$d^1_{0,0}\rightarrow d^2_{2,2}\rightarrow d^2_{1,2}\rightarrow d^2_{1,1}\rightarrow d^2_{0,2}\rightarrow d^2_{0,1}\rightarrow d^2_{0,0}$\\
\\
$d^2_{0,0}\rightarrow d^3_{3,3}\rightarrow d^3_{2,3}\rightarrow d^3_{2,2}\rightarrow d^3_{1,3}\rightarrow d^3_{1,2}\rightarrow d^3_{1,1}\rightarrow d^3_{0,3}\rightarrow d^3_{0,2}\rightarrow d^3_{0,1}\rightarrow
d^3_{0,0}$\\
$$\vdots$$\\
$d^{k-1}_{0,0}\rightarrow d^{k}_{k,k}\rightarrow d^k_{k-1,k}\to d^{k}_{k-1,k-1}\to\cdots\to d^{k}_{i,k}\to
d^k_{i,k-1}\to\cdots\to d^k_{i,i}\to\cdots\to d^k_{0,k}
\to\cdots\to d^k_{0,0}$ \\
$$\vdots$$

\medskip

\noindent The function $f$ is defined on the limit and isolated points  as follow:
\begin{enumerate}
\item $f(d^{k}_{i,l})=d^{k}_{i,l-1}$ for every $i\in\mathbb N$ and $k\geq l>i$.

\item $f(d^{k}_{i,i})=d^{k}_{i-1,k}$ for every $i\in\mathbb N$ and $k>i$.

\item $f(d_i) = d_i$ for every $i\in\nat$.

\item  $f(d_{i,l}) =  d_{i,l-1}$ for every $i\in\nat$ and $l>i$.

\item $f(d_{i,i}) =d_{i-1}$ for every $i>0$.

\item $f(d_{0,0})=d$.

\item $f(d)= d$.
\end{enumerate}
It is not difficult to verify that
\begin{enumerate}
  \item[$(a)$] $f[D_{i,l}]=D_{i,l-1} $, for every $i\in\nat$ and $l>i+1$, and $f[D_{i,i+1}]=D_{i,i}\setminus\{d^i_{i,i}\} $.
  \item[$(b)$] $f[D_{i,i}]=\{d^{k}_{i-1,k}:k>i\}$ for all $i>0$.
  \item[$(c)$] $f[D_{0,0}]=\{d^{k}_{k,k}:k\geq 1\}$.
\end{enumerate}

First, we will show that $f$ is continuous on $\omega^3+1$.
Clearly, by  identities $(6)$ and $(c)$, we obtain that $f$ is continuous on $d_{0,0}$.
Note that conditions $(4)$, $(5)$,
$(a)$ and  $(b)$ guarantee  that $f$ is continuous at $d_{i,l}$
for each $0<i\leq l$. Since $d_l$ is fixed for each
$l\in\mathbb N$, then  $(a)$ implies that $f$ is
continuous at $d_l$. It is not hard to verify the continuity of $f$ at $d$.

Finally, let us check that $f^p$ is discontinuous on $\omega^3+1$,
for each $p\in\nat^*$. Fix $p\in\nat^*$ and $1 \leq i$. We know that  $(d_{i,l})_{i<l}$
is a strictly increasing sequence converging to $d_i$. By condition $(4)$ and
$(5)$, for each $i<l$ there exist $n_l\in\nat$ such that
$f^m(d_{i,l})=d_{i-1}$ for $m\geq n_l$. Consequently,
$f^p(d_{i,l})=d_{i-1}$ for each $l<i$, but $d_i$ is fixed and so $f^p(d_i)=d_i$.
Therefore, $f^p$ is discontinuous at $d_i$. $\Box$

\medskip

We finish with a list of  open questions that the authors were not able to solve.

\medskip

The functions given in Examples  \ref{familia disjunta2}, \ref{orbita densa 2} and \ref{orbitad-densa3} have a dense orbit. It is then natural to ask:

\begin{question} Giving a countable ordinal $\alpha > 3$, is it possible to define a continuous function  $f: \omega^\alpha+1 \to \omega^\alpha+1$  with dense orbit?
\end{question}

The following question is related to  Theorem \ref{orbita-densa-omega2} and Example \ref{orbitad-densa3}.

\begin{question} Given a discrete dynamical system $(\omega^\alpha+1,f)$  with dense orbit, where $\alpha\geq 3$ is a countable  ordinal, is $E(\omega^\alpha+1,f)$ always countable?
\end{question}

It is not hard to see that $E(\omega^\alpha+1,f)$ could be a convergent sequence with its limits point without regarding the size of $\alpha \geq 1$. This leads to formulate the next question.

\begin{question} Let  $\alpha, \beta > 3$ be  countable  ordinals. Is there a continuous function $f: \omega^\alpha+1 \to \omega^\alpha+1$ such that $E(\omega^\alpha+1,f)$ is homeomorphic to
$\omega^\beta +1$?
\end{question}

The previous question may be stated in a more general form:

\begin{question} Given a compact metric countable  space $X$, is there a continuous function $f: X \to X$ such that
$E(X,f)$ is homeomorphic to $X$?
\end{question}

According to Theorem \ref{infipf}, if $E(X,f)$ is countable, then $P_f$ must be finite. Hence, we may ask:

\begin{question} Given a compact metric countable dynamical system $(X,f)$, if $P_f$ is finite, must $E(X,f)$ be countable?
\end{question}

Examples \ref{familia disjunta2} and \ref{orbita densa 2} provide two different functions $f_1,f_2:\omega^2+1\rightarrow\omega^2+1$ such that $E(\omega^2+1, f_i)$ is homeomorphic to $\omega^2+1$ for $i=1,2$.  The remarks following the proof of Theorem \ref{orbita-densa-omega2} assert that  $E(\omega^2+1,f_i)$, for $i = 1, 2$, is not algebraically isomorphic to $\omega^2 +1$ equipped with ordinal addition.   These facts naturally suggests the next question that could be easy to answer.

\begin{question}
Are the Ellis semigroups constructed in \ref{familia disjunta2} and \ref{orbita densa 2} algebraically homeomorphic?
\end{question}

We know, by Theorem \ref{infipf}, that  if $P_f$ is  infinite and $X$ is countable, then $E(X,f)$  is homeomorphic to $\cantor$. But what about the algebraic operation:

\begin{question}
Suppose $X$ is compact metric countable space and  $f,g: X \to X$ are continuous functions such that $P_f$ and $P_g$ are infinite.
Are the Ellis semigroups $E(X,f)$ and $E(X,g)$  algebraically homeomorphic?
\end{question}

\bibliographystyle{amsplain}


\end{document}